\documentclass{article}
\usepackage{graphicx} 
\usepackage{amssymb}
\usepackage{amsmath,amsthm, bbm, mathdots}
\usepackage{color, xcolor}
\usepackage{graphicx}
\usepackage{setspace}
\usepackage{xspace}
\usepackage{multirow}
\usepackage{array}
\usepackage{complexity}
\usepackage{geometry}
\usepackage{mathtools}
\usepackage{algorithm,algpseudocode}
\usepackage{dirtytalk}
\usepackage{tikz}
\usepackage
{hyperref}
\usepackage[capitalize]{cleveref}
\usepackage{dsfont}

\hypersetup{colorlinks=true,urlcolor=blue,linkcolor=magenta,citecolor=[rgb]{.42,.56,.14},}

\usepackage{indentfirst}
\makeatletter
\renewcommand{\@seccntformat}[1]{\csname the#1\endcsname.\quad}
\makeatother

\theoremstyle{plain}
\newtheorem{theorem}{Theorem}

\newtheorem{lemma}[theorem]{Lemma}
\newtheorem{claim}[theorem]{Claim}

\newtheorem{definition}[theorem]{Definition}

\theoremstyle{remark}

\theoremstyle{plain}
\newclass{\DNF}{DNF}
\newclass{\DNFs}{DNFs}
\newclass{\ACzero}{AC^0}
\newclass{\TCzero}{TC^0}

\renewcommand{\Pr}{\mathop{\bf Pr\/}}

\renewcommand{\epsilon}{\varepsilon}

\usepackage{amssymb,amsmath,csquotes}

\title{The Story of Sunflowers}

\author{Anup Rao \\ University of Washington \\ \texttt{anuprao@cs.washington.edu}}

\begin{document}

\maketitle

\begin{abstract}
Sunflowers, or $\Delta$-systems, are a fundamental concept in combinatorics introduced by  Erd\H{o}s and Rado in their paper: {\em Intersection theorems for systems of sets}, J. Lond. Math. Soc. (1) {\bf 35} (1960), 85--90. A sunflower is a collection of sets where all pairs have the same intersection.
This paper explores the wide-ranging applications of sunflowers in computer science and combinatorics. We discuss recent progress towards the sunflower conjecture and present a short elementary proof of the best known bounds for the robust sunflower lemma.
\end{abstract}

\section{Introduction}
Elegant definitions captivate mathematicians, and their appeal deepens when they lead to unexpected applications. Such is the story of \emph{sunflowers}, first studied by Erd\H{o}s and Rado in their paper ``Intersection theorems for systems of sets" in 1960 \cite{ErdosRado1960}. The term sunflower was later  coined by Frankl and  first appeared in the work of Deza and Frankl \cite{DezaFrankl1981}; Erd\H{o}s and Rado originally used the term ``$\Delta$-system" for the same concept.
\begin{figure}[h]
    \centering \includegraphics[width=0.3\textwidth]{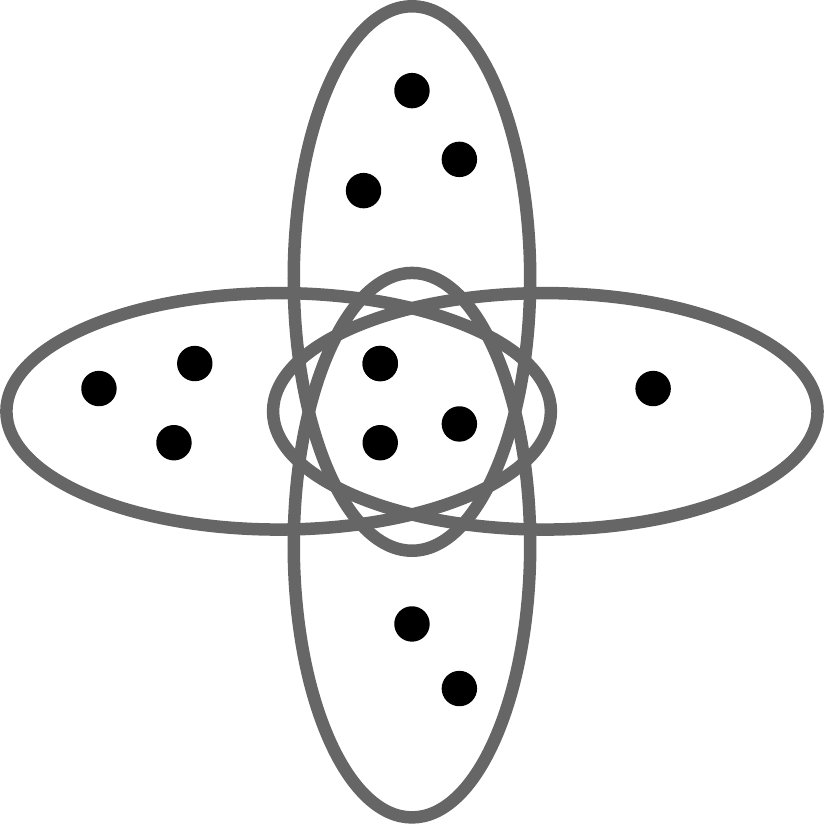}
    \caption{A sunflower with $4$ petals.}
    \label{fig:sunflower}
\end{figure}

A sunflower with $w$ petals is a collection of $w$ sets such that every pair of sets has the same intersection. The common intersection is called the \emph{core} or the \emph{kernel}. Erd\H{o}s and Rado proved that every large collection of sets must contain a sunflower. They gave a simple inductive argument showing that every collection of more than $k! \cdot (w-1)^k$ sets of size at most $k$ contains a sunflower with $w$ petals. The proof is nice:  consider a collection of $k! \cdot (w-1)^k$ sets of size $k$. Select a set $S_1$, and then a set $S_2$ disjoint from $S_1$, and then $S_3$ disjoint from $S_1,S_2$, and continue in this way. If we end up with $w$ disjoint sets, they form the desired sunflower. If this process fails, it must be because we have found $t<w$ sets $S_1,\ldots,S_t$ such that every set in the family intersects one of the sets $S_i$. But then there must be some element of $S_1 \cup \dotsb \cup S_t$ that is contained in at least $\frac{k! \cdot (w-1)^k}{|S_1 \cup \dotsb \cup S_t|} \geq (k-1)! \cdot (w-1)^{k-1}$ of the sets in the family. This family of  $(k-1)! \cdot (w-1)^{k-1}$ sets is effectively a family of sets of size $k-1$, because they all share an element. So, we repeat the argument to find a sunflower in this new family, adding the common element to the core.

Is this bound tight? Erd\H{o}s and Rado included a variant of the following example, which is a family of  $(w-1)^k$ sets of size $k$ that does not have a   sunflower with $w$ petals. Let $[n]=\{1,\ldots,n\}$, and consider all subsets $S \subseteq [(w-1)k]$ such that, for each $a = 1, \ldots, k$, the set $S$ contains exactly one element $s$ satisfying $s \equiv a \pmod{k}$. Consider $w$ distinct sets forming a supposed sunflower. For any $a$, there are only $w-1$ elements congruent to $a \pmod{k}$. By the pigeonhole principle, two of the $w$ petals must share an element congruent to $a \pmod{k}$, which must belong to all petals by the sunflower property.  Since this holds for all $a$, the $w$ sets must be identical, contradicting their distinctness. Hence, no sunflower with $w$ petals exists in this family.  Erd\H{o}s and Rado conjectured that any family of sets of size at most $k$ avoiding a $w$-petal sunflower has at most $C^k$ sets, where $C$ is a constant depending only on $w$ and independent of $k$.

The study of sunflowers is part of the broader saga of Ramsey theory, whose moral is that surprising structures can be found in large objects.  Notable early examples include Van der Waerden’s theorem (1927) \cite{vanDerWaerden1927}, which guarantees a monochromatic arithmetic progression in any coloring of $[n]$ for large $n$, Sperner’s theorem (1928) \cite{Sperner1928}, which ensures that a large family of subsets of $[n]$ contains a pair where one is a subset of the other, and Ramsey’s theorem (1929) \cite{Ramsey1929}, which states that a graph on $n$ vertices has either a clique or independent set of size $(1/2)\log_2 n$. Other landmarks include Tur\'{a}n’s theorem (1941) \cite{Turan1941}, Roth’s theorem (1953) \cite{Roth1953}, the Erd\H{o}s–Ko–Rado theorem (1961) \cite{ErdosKoRado1961}, the Hales–Jewett theorem (1963) \cite{HalesJewett1963}, Ajtai and Szemer\'{e}di’s corner theorem (1974) \cite{AjtaiSzemeredi1974}, and Szemer\'{e}di’s theorem (1975) \cite{Szemeredi1975}. Each of these results uncovers inevitable patterns in large structures, so it was natural for Erd\H{o}s and Rado to look for sunflowers in 1960. Yet, they could scarcely have foreseen  the immense number of applications of sunflowers, particularly in computer science, that would arise decades later. We discuss two applications in \Cref{apps}. 

\subsection{Progress}
Erd\H{o}s and Rado's conjecture remains open, despite extensive efforts to resolve it. In 1997, Kostochka \cite{Kostochka1997} proved  that a family of $O(k! \cdot (\frac{\log \log \log k}{\log \log k})^k)$ sets of size $k$  contains a sunflower with $3$ petals, showing that Erd\H{o}s and Rado's bound can be improved. In 2019, Alweiss, Lovett, Wu and Zhang \cite{alweiss2021} made a significant breakthrough. They proved that a family of $O(w^3 \log k \log \log k)^k$ sets of size $k$ must contain a sunflower with $w$ petals. As is often the case, progress came through finding a stronger structure. In this case, the stronger structure is called a \emph{robust sunflower}, a concept first studied by Rossman \cite{Rossman2014}. Let $A(\gamma)$ denote a random subset sampled by including each element of the universe in $A(\gamma)$ independently with probability $\gamma$.
\begin{definition}
    A family $\mathcal{R}$ of at least $2$ subsets of $[n]$ is a $(\gamma,\epsilon)$-\emph{robust sunflower} if, letting $C$ denote the intersection of all sets in $\mathcal{R}$, we have $$\Pr[\exists S \in \mathcal{R}, S \subseteq A(\gamma) \cup C] \geq 1-\epsilon.$$ 
\end{definition}
The set $C$ in the definition is the analogue of the core of the sunflower, and we call it the core in the context of robust sunflowers as well. However, unlike with sunflowers, two sets of the robust sunflower may intersect outside the core. The breakthrough of \cite{alweiss2021}, gave a new robust sunflower lemma, and subsequent work by myself \cite{Rao20}, Frankston, Kahn, Narayanan, and Park \cite{FrankstonKNP19}, and Bell, Chueluecha, and Warnke \cite{BellCW21} (see the blog post by Tao \cite{Tao}) gave the following version:
\begin{lemma}[Robust Sunflower Lemma]\label{robust}
There is a universal constant $c > 1$ such that if $\epsilon \leq 1/2$, then every family of $( c\log ( k/\epsilon )/\gamma)^k$ sets of size at most $k$ contains a $(\gamma,\epsilon)$-robust sunflower.
\end{lemma}
Every $(1/(2w),1/2)$-robust sunflower contains a standard sunflower with $w$ petals \cite{BellCW21}. Indeed, if we  partition the universe uniformly randomly into $2w$ disjoint sets $B_1,\ldots,B_{2w}$, the robust sunflower lemma guarantees that we expect to find at least  $(2w)/2=w$ sets $S_1,\ldots,S_w \in \mathcal{R}$ and $w$ disjoint parts $B_{j_1},\ldots,B_{j_w}$ with $C \subsetneq S_i \subseteq B_{j_i} \cup C$. The sets $S_1,\ldots,S_w$ must form a sunflower with $w$ petals and core $C$. So, \Cref{robust} implies that  every family of $O(w \log k)^k$ sets of size $k$ must contain a sunflower with $w$ petals. The only difference between this bound and the  conjecture of Erd\H{o}s and Rado is the presence of the $\log k$ term. This dependence is necessary for robust sunflowers, as shown by  \cite{alweiss2021}, who found  a family of $\Omega(\log k)^{k-\sqrt{k}}$ sets that does not contain a $(1/2,1/2)$-robust sunflower. Nevertheless, it is quite possible that the sunflower conjecture of Erd\H{o}s and Rado holds in its original form.

The proofs of \Cref{robust} given in  \cite{Rao20, BellCW21, Tao} relied on ideas from information theory. We give an elementary proof of \Cref{robust} in  \Cref{sectionrobust}, which improves on an elementary proof of the classical sunflower lemma given in my earlier exposition \cite{Rao2023}.

\section{Applications of sunflowers} \label{apps}
Sunflowers reveal patterns that prove valuable in many computer science applications. Here, we outline two interesting applications.

\subsection{Dynamic data structures}
Our first application demonstrates how sunflowers provide lower bounds for dynamic data structures~\cite{RamamoorthyRao2018,FrandsenMS1997,GalM2007}. Consider an algorithm tasked with maintaining a subset $T \subseteq [m]=\{1,\ldots,m\}$ and computing its minimum element. The algorithm uses $n$ memory cells, each storing a single element from $[m]$. It must support adding or deleting elements from $T$ and computing the minimum element of $T$, and each operation should involve at most $t$ accesses to memory cells. The operation of adding or deleting an element $i$ from $T$ is performed by accessing a predetermined set of memory cells $S_i \subseteq [n]$ with $|S_i|\leq t$. More precisely, when element $i$ is added to or removed from $T$, the algorithm may only read from and write to the memory cells corresponding to  $S_i$. The goal is to minimize $t$ to improve efficiency. A well-known data structure called a binary search tree achieves $t \lesssim \log m$. Using the sunflower lemma, Ramamoorthy and I showed that every data-structure must satisfy $t \gtrsim \log m / \log \log m$~\cite{RamamoorthyRao2018}. It is an open problem to understand what can be achieved by adaptive data structures that are not limited to accessing a fixed set of memory cells $S_i$. 

The connection to sunflowers is easiest to explain in the special case where the minimum is computed by accessing a fixed set of memory cells $S_0 \subseteq [n]$. The main observation is that if the sets $S_1 \cup S_0, \ldots, S_w \cup S_0$ happen to form a sunflower with a common intersection $C$, the algorithm is severely constrained about how it stores  information about subsets $T \subseteq [w]$. Consider storing the set $[w]$ and then deleting the elements of the set $[w] \setminus T$, so that the memory now encodes $T$. After these operations, we claim that the contents of the memory cells in the common intersection $C$ are sufficient to recover $T$. Indeed, $T$ can be recovered by repeatedly computing and deleting the minimum of $T$ using only the memory cells in $C$. The sunflower structure ensures that for $i \in T$, the contents of $S_i\setminus C$ are not modified after the set $[w]$ is added using the algorithm, so $T$ is encoded by the memory cells of $C$. This means that the size of the common intersection $|C|$ must be large enough to encode all $2^w$ subsets of $[w]$, in other words, $m^{|C|}\geq 2^w$, which implies that $2t \geq |C| \geq w / \log m$. The sunflower lemma ensures that if $m$ is sufficiently large and $t$ is small, a sunflower with $w$ petals exists among the sets $S_1 \cup S_0,\ldots, S_m \cup S_0$. A straightforward calculation using the bounds given by the sunflower lemma yields $t \geq \Omega(\log m / \log \log m)$.

\subsection{Monotone circuit complexity}
Our second application is to  proving exponential lower bounds for computing natural functions using monotone circuits, following classical arguments due to Razborov \cite{Razborov1985}, and Alon and Boppana \cite{AlonBoppana1987} and the very recent work of Cavalar, G{\"o}{\"o}s, Riazanov, Sofronova, and Sokolov  \cite{Cavalar2025}. A monotone circuit is a computational model represented as a directed acyclic graph (DAG). The input nodes of the DAG correspond to Boolean variables $x_1, \ldots, x_n$, each taking a value in $\{0, 1\}$ representing specific information, such as the presence or absence of an edge in a graph. Each internal node computes either an AND (logical conjunction) or OR (logical disjunction) operation on the outputs of two previous nodes, and the circuit has a single output node that produces the final result. Unlike general Boolean circuits, monotone circuits are restricted to AND and OR operations, with no negations allowed, meaning they can only compute monotone Boolean functions---functions where changing an input from $0$ to $1$ never flips the output from $1$ to $0$. The circuit's size is defined as the number of gates (AND or OR nodes). Proving similar exponential lower bounds for general Boolean circuits (which allow negations) would resolve the $\mathsf{P}$ vs $\mathsf{NP}$ problem,  a major open question in computer science. These results suggest that the use of negations in algorithms is a major barrier to proving $\mathsf{P} \neq \mathsf{NP}$.

Here, we describe lower bounds for computing two natural graph properties using monotone circuits. The inputs $x_1, \ldots, x_n$ encode a graph on $m$ vertices, with $n = \binom{m}{2}$ and  $x_i$ indicates whether the $i$-th edge is present. In the first example, the goal is to determine whether the graph contains a \emph{clique} of size $k\approx n^{1/3}$: a set of  $k$ vertices fully connected by edges. In the second example, the goal is to determine whether a  graph with an even number of vertices $m$ contains a \emph{perfect matching}, namely a set of $m/2$ disjoint edges. We give rough sketches showing how to prove that any monotone circuit computing these functions requires exponentially many gates.

\subsection{Cliques} Here we present the approach of \cite{Razborov1985, AlonBoppana1987}. 
Suppose a small monotone circuit exists computing whether the graph has a clique of size $k$. Let $G$ be a random graph containing a clique of size $k$ and no other edges, meaning only the edges within the clique are present. Let $H$ be a random $(k-1)$-partite graph, where vertices are partitioned into $k-1$ parts, and only the edges connecting distinct parts are present. This ensures $H$ never contains a clique of size $k$. Finally, let $F$ be a random graph which is equal to $G$ with probability $1/2$, and equal to $H$ with probability $1/2$.

For a set of vertices $S$, let $\mathsf{clique}_S(x)$ denote the Boolean function that outputs 1 when $x$ contains a clique on $S$. The key claim of the proof is the following, which we state here without being too specific about the parameters:
\begin{claim}\label{cliqueclaim}
   If a small monotone circuit $f$ computes whether a graph has a clique of size $k$, then there is a small collection of sets $S_1,\ldots,S_T$ such that $$f(F) = \bigvee_{i=1}^T \mathsf{clique}_{S_i}(F)$$ holds with high probability over the choice of $F$.
\end{claim}
It is not hard to see that any expression that is an OR of clique functions of this type must have $T$ that is exponentially large if it computes the clique function on $F$ correctly. So, the heart of the matter is finding a small formula as above that approximates the computation carried out by the monotone circuit. To carry out the approximation, we proceed by induction. For each gate $g$ computed by the circuit, we find a collection of  sets $S_1, \ldots, S_T \subseteq [m]$ such that,
\begin{align}
g(G) \leq  \bigvee_{j=1}^T \mathsf{clique}_{S_j}(G), \label{fapprox1}
\end{align}
and
\begin{align}
g(H) \geq  \bigvee_{j=1}^T \mathsf{clique}_{S_j}(H) \label{fapprox2}
\end{align}
both hold with high probability over the choice of $G,H$. Such an approximation holds for input variables because each $x_i$ is equal to $\mathsf{clique}_e$, where $e$ is the $i$'th edge of the graph. If the circuit takes the OR of two expressions of size $T_1,T_2$ as in \Cref{fapprox1,fapprox2}, the result can be computed by a similar expression with $T = T_1+T_2$. If the circuit takes the AND of such expressions, the result can be computed by a similar expression of size $T=T_1 \cdot T_2$ using the distributivity of AND and OR, because the choice of $G,H$ implies that
\begin{align*}
\mathsf{clique}_{S_1}(G) \wedge \mathsf{clique}_{S_2}(G) = \mathsf{clique}_{S_1 \cup S_2}(G), 
\end{align*}
and
\begin{align*}
\mathsf{clique}_{S_1}(H) \wedge \mathsf{clique}_{S_2}(H) \geq  \mathsf{clique}_{S_1 \cup S_2}(H), 
\end{align*}

The problem is that these ideas only give $T$ in \Cref{fapprox1,fapprox2} that is exponentially larger than the size of the monotone circuit. We use sunflowers to give efficient approximations. If any set $S_i$ in \Cref{fapprox1,fapprox2} is larger than $k$, we can simply remove $S_i$ from the collection. This only causes a problem when $S_i$ is contained in the clique of $G$, and that occurs with extremely small probability. So, we can assume that all sets $S_i$ are of size at most $k$. The key observation using sunflowers is as follows. If $w$ is large, and $S_1, \ldots, S_w$ forms a sunflower with core  $C$, then
\[
\bigvee_{j=1}^w \mathsf{clique}_{S_j}(G) \leq \mathsf{clique}_C(G),
\]
by definition. Moreover, when the input is $H$,
\[
\bigvee_{j=1}^w \mathsf{clique}_{S_j}(H) \geq \mathsf{clique}_C(H)
\]
holds with high probability over the choice of $H$. That is because a clique on $C$ in $H$ likely implies a clique on at least one of the petals of the sunflower, since if the vertices of $C$ form a clique, the events that the petals are also cliques are independent, and the parameters are chosen such that each event occurs with significant probability. 
So, if $T$ is large, the OR of cliques in \Cref{fapprox1,fapprox2} can be simplified by using the sunflower lemma to find a sunflower, and replacing the sunflower with its core $C$, preserving \Cref{fapprox1,fapprox2} except with some small error. Repeatedly applying this process in a small circuit yields \Cref{cliqueclaim}.

\subsection{Perfect matchings}
The lower bound for computing whether the graph has a perfect matching has a similar high-level structure to the proof for cliques \cite{Cavalar2025}.  Let $G$ be a graph containing exactly one uniformly random perfect matching. Let $H$ be a graph that is the union of two uniformly random disjoint odd-sized cliques covering all vertices.
Thus, $G$ always has a perfect matching, while $H$ does not, because an odd-sized clique cannot have a perfect matching. We set $F$ to be equal to $G$ with probability $1/2$, and equal to $H$ with probability $1/2$.

We approximate the circuit’s computation using indicators for matchings. Let $M$ denote a matching, namely a set of pairwise disjoint edges in the graph, and let $1_{M}$ denote the function that outputs $1$ only when the matching $M$ is present. The key claim is:

\begin{claim}\label{matchingclaim}
   If a small monotone circuit $f$ computes the perfect matching function, there exists a small set of matchings $\{M_1,\ldots,M_T\}$ such that $$f(F) = \bigvee_{i=1}^T 1_{M_i}(F)$$ holds with high probability over the choice of $F$.
\end{claim}
It is not hard to verify that a small OR of matchings cannot possibly compute the perfect matching function correctly on $F$. For each gate $g$ of the circuit, we construct a set of matchings $\{M_1, \ldots, M_T\}$ such that  
\begin{align}
    g(G) \leq \bigvee_{i=1}^T 1_{M_i}(G) \label{mapprox1}
\end{align}
and
\begin{align}
    g(H) \geq  \bigvee_{i=1}^T 1_{M_i}(H) \label{mapprox2}
\end{align} 
both hold with high probability over the choice of $G,H$.
The input variables correspond to matchings of size $1$, and we can generate \Cref{mapprox1,mapprox2} inductively for each gate  using the distributivity of AND and OR, because the choice of $G,H$ means 
\begin{align*}
    1_{M_1}(G) \wedge 1_{M_2}(G) = \begin{cases}
        0 & \text{if $M_1 \cup M_2$ is not a matching}\\
        1_{M_1 \cup M_2}(G) &\text{otherwise.}
    \end{cases}
\end{align*}And,
\begin{align*}
    1_{M_1}(H) \wedge 1_{M_2}(H) \geq  \begin{cases}
        0 & \text{if $M_1 \cup M_2$ is not a matching}\\
        1_{M_1 \cup M_2}(H) &\text{otherwise.}
    \end{cases}
\end{align*}
Thus, we have 
\begin{align*}
\Big(\bigvee_{i=1}^{T_1} 1_{M_i}(G)\Big)\wedge \Big(\bigvee_{j=1}^{T_2} 1_{M_j}(G)\Big) = \bigvee_{i,j: M_i \cup M_j (G) \text{ is a matching}} 1_{M_i \cup M_j}(G),
\end{align*}and 
\begin{align*}
\Big(\bigvee_{i=1}^{T_1} 1_{M_i}(H)\Big)\wedge \Big(\bigvee_{j=1}^{T_2} 1_{M_j}(H)\Big) \geq  \bigvee_{i,j: M_i \cup M_j(H) \text{ is a matching}} 1_{M_i \cup M_j}(H).
\end{align*}

The problem is that the value of $T$ obtained in \Cref{mapprox1,mapprox2}  this way is too large to give a result. We use sunflowers to give efficient approximations to such expressions. If any matching $M_i$ occurring in \Cref{mapprox1,mapprox2} involves more than $k$ edges, we can simply remove it from the expression. This only causes a problem when the random matching of $G$ contains $M_i$, and that happens with very low probability. Repeat this until all matchings are of size at most $k$.

Matchings do not directly correspond to vertex sets. A key insight is to apply the robust sunflower lemma by
first identifying a large \emph{blocky} subset of the matchings given in \Cref{mapprox1,mapprox2}. A set of matchings, is blocky if there exist $k$ disjoint cliques such that each matching’s edges lie in distinct cliques. Given $T$ matchings of size at most  $k$, some subset of roughly $(1/k)^{2k} \cdot T$ matchings must be blocky. That is because if we sample a random partition of the vertices into $k$ parts, then for any matching of size $k$, the first edge will be contained in the first part with probability $1/k^2$, the second edge  will be contained in the second part with probability $1/k^2$,  and so on. The entire matching will be consistent with the random partition with probability $(1/k)^{2k}$, so by averaging, there must be one partition that is consistent with $(1/k)^{2k} \cdot T$ matchings. 

There is a one-to-one correspondence between each matchings of this blocky set and the set of $2k$ vertices that it touches. If we associate each matching $M_i$ with the set of vertices $S_i$ that it touches, and if $T$ is large enough in \Cref{mapprox1,mapprox2}, we can apply the robust sunflower lemma, \Cref{robust}, to the corresponding family of subsets of the vertices to find a robust sunflower $\mathcal{R}$ corresponding to a subset of the matchings. Let $C$ denote the core of the robust sunflower. Let $M$ be the unique matching that is consistent with the blocky partition whose edges are completely contained in $C$. Then 

\begin{align*}
    \bigvee_{i \in \mathcal{R}} 1_{M_i}(G) \leq 1_M(G) 
\end{align*}
by construction, and one can use the robust sunflower property to prove 
\begin{align*}
    \bigvee_{i \in \mathcal{R}} 1_{M_i}(H) \geq 1_M(H) 
\end{align*}
holds with high probability. As before, this gives a subroutine to find efficient approximations. By replacing
\begin{align*}
    \bigvee_{i \in \mathcal{R}} 1_{M_i}(F) = 1_M(F), 
\end{align*} we obtain a smaller set of matchings for \Cref{mapprox1,mapprox2}. Repeating this ultimately yields \Cref{matchingclaim}.

\subsection{Other applications}
Beyond these, sunflowers have revealed structure in large systems. For example, Erd\H{o}s and S\'{a}rk\"{o}zy used sunflowers to prove that every large set of integers contains subsets whose sums form an arithmetic progression \cite{ErdosS}. The concepts surrounding sunflowers and recent progress in \cite{alweiss2021} also resolved conjectures about the \emph{threshold of monotone functions}. For a monotone function $f:\{0,1\}^n \to \{0,1\}$, the threshold is the probability $p$ such that if each input is independently set to 1 with probability $p$, then $f$ equals 1 with probability $1/2$. This parameter is significant in studying random processes. Talagrand \cite{Talagrand} and Kahn and Kalai \cite{KahnKalai2007} proposed methods to estimate the threshold of monotone functions and made conjectures about these estimates. These were confirmed by Frankston, Kahn, Narayanan, and Park \cite{FrankstonKNP19}, and Park and Pham \cite{ParkPham2022}, building on techniques from \cite{alweiss2021}. Further details can be found in \cite{Rao2023}.

\section{The robust sunflower lemma} \label{sectionrobust}
Recall the statement of the robust sunflower lemma:
{
\renewcommand{\thetheorem}{\ref{robust}}
\begin{lemma}[Robust Sunflower Lemma]
There is a universal constant $c > 1$ such that if $\epsilon \leq 1/2$, then every family of $( c\log ( k/\epsilon )/\gamma)^k$ sets of size at most $k$ contains a $(\gamma,\epsilon)$-robust sunflower.
\end{lemma}
\addtocounter{theorem}{-1}
}

In this section, we present an elementary proof of the robust sunflower lemma. 
Let $r = c\log ( k/\epsilon )/\gamma$, where the logarithm is computed in base $2$, and $c$ is to be determined. Let  $\mathcal{F}$ be a family of $r^k$ distinct sets. We prove by induction on $k$ that $\mathcal{F}$ contains a $(\gamma,\epsilon)$-robust sunflower.

When $k=1$, we prove that every family of $r\geq 2$ sets is a robust sunflower with core $C=\emptyset$. If the family contains the empty set, then the conclusion trivially holds. Otherwise, the probability that the random set $A(\gamma)$ does not include one of the $r$ sets of the family is at most $(1-\gamma)^r\leq e^{-\gamma r} \leq e^{-c \log(1/\epsilon)}\leq \epsilon$
for $c$ large enough. Here we used the identity $1+x \leq e^{x}$.

For $k>1$, if there exists a set $Z$ with $0< |Z|<k$ contained in some $r^{k-|Z|}$ sets of $\mathcal{F}$, we apply induction on the family of sets containing $Z$, which are effectively sets of size at most $k-|Z|$. We add the elements of $Z$ to the core $C$ found by induction to obtain our final core. Otherwise, it must be that for every non-empty set $Z$, the number of sets of the family containing $Z$ is at most $r^{k-|Z|}$. In particular, the intersection of all the sets in the family must be empty. We prove that such a family is a $(\gamma,\epsilon)$-robust sunflower with core  $C=\emptyset$. The family certainly has $2$ distinct sets since $r^k\geq 2$. All that remains is to prove $\Pr[\exists S \in \mathcal{R}, S \subseteq A(\gamma)] \geq 1-\epsilon$.

Set $\ell = \lceil c_0 \cdot  \log(k/\epsilon)\rceil$, for a  constant $c_0>0$ to be determined. We add dummy elements to the universe until $n$ is large enough so that $\Pr[|A(\gamma)| < \gamma n/2] < \epsilon/2$. By the Chernoff bound, $\Pr[|A(\gamma)| < \gamma n/2] \leq e^{-\Omega(\gamma n)}$, so there is a value of $n$ for which this is achieved. If $|A(\gamma)|\geq  \gamma n/2$, sample a uniformly random set $B \subseteq A(\gamma)$ of size $\ell \cdot \lfloor\gamma n/(2\ell) \rfloor \leq \gamma n/2$. Conditioned on the event that $|A(\gamma)|\geq \gamma n/2$, $B$ has the distribution of a uniformly random set of size $\ell \cdot \lfloor\gamma n/(2\ell)\rfloor$. It is enough to prove:
\begin{align}\label{Bbound}
\Pr[\exists S\in  \mathcal{F}, S \subseteq B]\geq 1-\epsilon/2,
\end{align} because this implies that $$\Pr[\exists S\in  \mathcal{F}, S \subseteq A(\gamma)]\geq 1-\epsilon,$$ as required. To prove \Cref{Bbound}, partition  $B$ into $\ell$ sets uniformly at random: $B_1 \cup \dots\cup B_{\ell}= B$, with  $|B_i|= \lfloor\gamma n/(2\ell) \rfloor$. We denote $B^i = B_1  \cup \dots \cup B_i$, and $B^0=
\emptyset$.

\begin{figure}[h]
    \centering \includegraphics[width=0.5\textwidth]{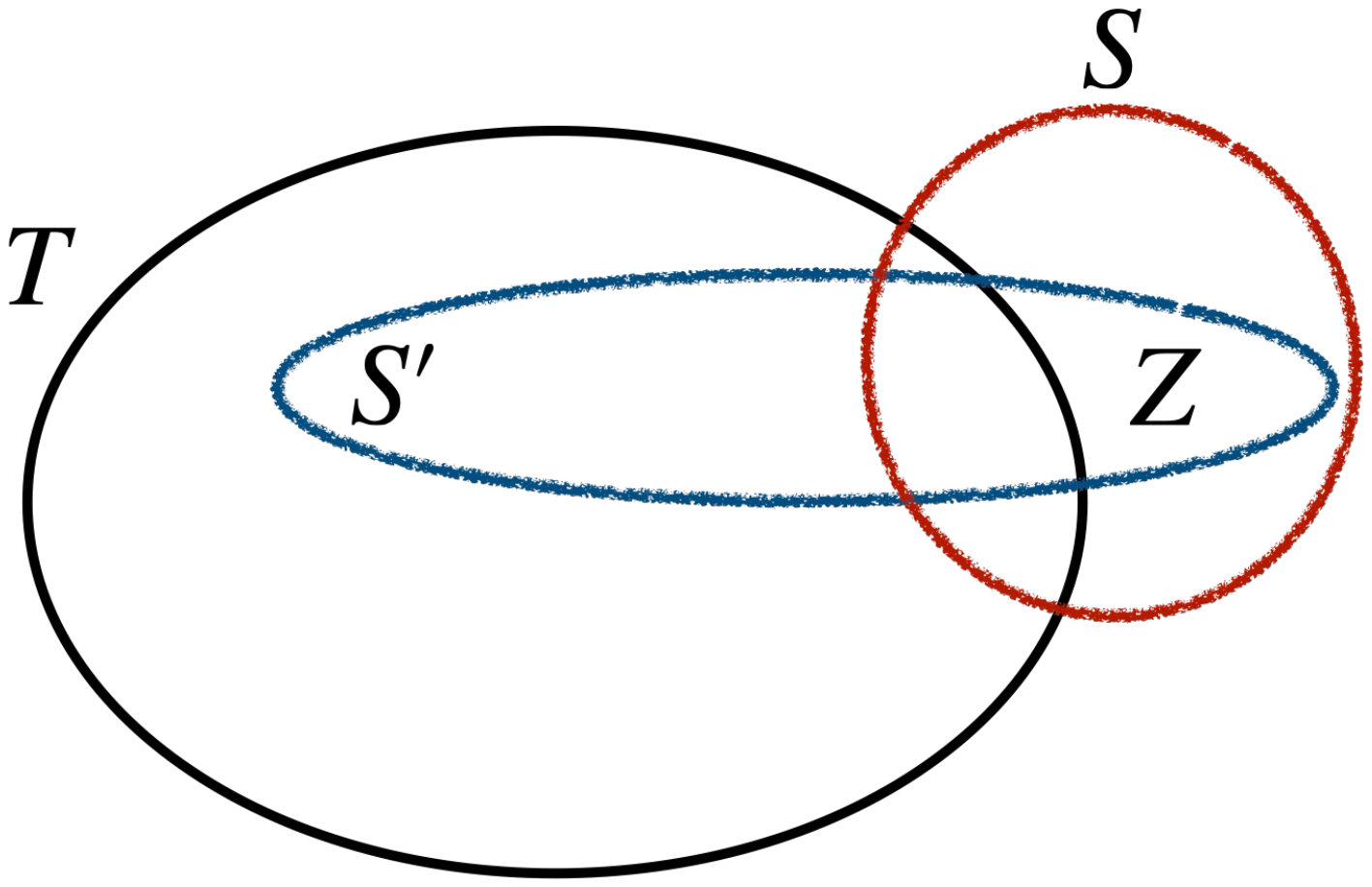}
    \caption{$S$ is $m$-covered by $T$ if $|Z|\leq m$.}
    \label{fig:mcover}
\end{figure}
Say that a set  $S \in \mathcal{F}$ is $m$-covered by a set $T \subseteq [n]$ if there is $S' \in \mathcal{F}$ such that $S' \subseteq T \cup S$ and $|S'\setminus T|\leq  m$. Observe that if a set is  $m$-covered by $T$, then it is also $m$-covered by all $U \supseteq T$. If a set is $0$-covered by $A$, then there must be    $S' \in \mathcal{F}$ with $S' \subseteq A$; our goal is to prove that a set is $0$-covered with probability at least $1-\epsilon$. Here is the key technical claim:
\begin{claim}\label{mainclaim}
    If $N$   sets of $\mathcal{F}$ are $m$-covered by $B^i$, then with probability at least $1/2$ over the choice of $B_{i+1}$, at least $N-  r^k/4^{m}$ sets are $(m/2)$-covered by $B^{i+1}$. 
\end{claim}

To prove the robust sunflower lemma, we repeatedly apply \Cref{mainclaim}. There are $r^k$ sets in $\mathcal{F}$, and each of them is $k$-covered by $B^0=\emptyset$. The claim guarantees that with probability $1/2$, at least $r^k-r^k/4^k$ sets are  $(k/2)$-covered by $B^1$. If this does not happen, it does not matter so much, because the  sets are still $k$-covered by $B^1$. In the next step, we can potentially find $r^k- r^k/4^k - r^k / 4^{k/2}$ sets that are $(k/4)$-covered, and so on. If $1+\lceil \log k \rceil$ rounds of this random process succeed, because  $k/2^{1+\lceil \log k \rceil}<1$, we are left with sets that are $0$-covered. 

There are $\ell$ rounds of sampling $B_1,\ldots,B_\ell$, and each succeeds with probability $1/2$. We have $1+\lceil \log k \rceil \leq \ell/3$, for $c_0$ large enough. By standard concentration bounds like Azuma's inequality, the probability that less than $\ell/3$ rounds succeed is at most $e^{-\Omega(\ell)} \leq \epsilon/2$, for $c_0$ chosen large enough. So, with probability at least $1-\epsilon/2$, we get $1+\lceil \log k \rceil$ successful steps. 

\Cref{mainclaim} guarantees that the number of sets that are $0$-covered at the end is at least (using the identity  $1/(1-x) = 1+x+x^2+\ldots$):
\begin{align*}
r^k -r^k(1/4^k + 1/4^{\lfloor k/2 \rfloor} +\ldots + 1 / 4) &\geq r^k- r^k ((1/4) + (1/4)^2 + \ldots) \\&= r^k-r^k (1/4)/(1-1/4)\\
&=  2r^k/3\\
&\geq 1.
\end{align*} So, the claim proves that $B^\ell$ contains some set of $\mathcal{F}$,  with probability at least $1-\epsilon/2$, establishing \Cref{Bbound}. It only remains to prove \Cref{mainclaim}. Fix $B^i$, and define the random family:
\[
\mathcal{Z} = \{ S \setminus B^{i+1} : \text{$S \in \mathcal{F}$ is $m$-covered by $B^i$, but not $(m/2)$-covered by $B^{i+1}$}\}.
\]
A set $Z \in \mathcal{Z}$ is \emph{minimal} if it does not contain any other set of  $\mathcal{Z}$. There cannot be a minimal set $Z$ with $|Z| \leq  m/2$, since the corresponding set $S$ is $(m/2)$-covered by $B^{i+1}$. We bound the expected number of minimal sets  of size $b \geq m/2$ using the following count: 

\begin{enumerate}
    \item Let $n' = n-|B^i|$, $t = |B_{i+1}|=\lfloor \gamma n /\ell \rfloor$. There are at most $\binom{n'}{t + b}$ choices for $Z \cup B_{i+1}$ with $|Z| = b$. We have:
    \[
    \binom{n'}{t + b} = \binom{n'}{t} \cdot \prod_{j=1}^b \frac{n' - t - j}{t + j} \leq \binom{n'}{t} \cdot \left( \frac{n}{t} \right)^b \leq \binom{n'}{t} \cdot (r/1000)^b,
    \]
since $\frac{n}{t} = \frac{n}{\lfloor \gamma n /\ell \rfloor}\leq \frac{2n}{\gamma n /\ell} \leq \frac{2c_0}{\gamma} \cdot \lceil \log (k/\epsilon)\rceil \leq r/1000$, if $c/c_0$ is large enough. 
    \item Given $Z \cup B_{i+1}$, let $S \setminus B^{i}$ be the smallest set (breaking ties by picking the lexicographically first one) among all $S \in \mathcal{F}$ that are $m$-covered by $B^i$ and satisfy $S \setminus B^{i} \subseteq Z \cup B_{i+1}$. This choice ensures that $|S \setminus B^i| \leq m$. We must have $Z \subseteq S \setminus B^{i}$, since otherwise $S \setminus B^{i+1}$ would be a proper subset of $Z$, contradicting $Z$'s minimality. Thus, there are at most $2^{m}$ choices for $Z$ given $Z \cup B_{i+1}$.
\end{enumerate}

Because there are exactly $\binom{n'}{t}$ choices for $B_{i+1}$, the above count shows that the  expected number of minimal sets of size $b$ is at most $2^m \cdot (r/1000)^b$. By assumption, each such minimal $Z$ of size $b$ is contained in at most $r^{k - b}$ sets of $\mathcal{F}$. Recall that there are no minimal sets of size less than $m/2$ in $\mathcal{Z}$. Every set that is not $(m/2)$-covered contains such a minimal $Z$. Thus, the expected number of sets that are not $(m/2)$-covered by $B^{i+1}$ is at most 
\begin{align*}
\sum_{b=\lceil m/2 \rceil}^\infty r^{k - b} \cdot 2^{m} \cdot (r/1000)^b &\leq r^k \cdot 2^m \cdot (1/1000)^{\lceil m/2 \rceil}\cdot  \sum_{j=0}^\infty (1/1000)^j\\
&= r^k \cdot 2^m \cdot (1/1000)^{\lceil m/2 \rceil}\cdot  (1/(1-1/1000))\\
&= r^k \cdot (2/\sqrt{1000})^m \cdot  (1000/999)\\
&\leq r^k/8^m.
\end{align*}
By Markov's inequality, we conclude that the probability that $r^k/4^m \geq 2 \cdot r^k/8^m$ sets are not $(m/2)$-covered is at most $1/2$, proving \Cref{mainclaim}.

\section{Acknowledgements}
Thanks to Paul Beame, Bruno Cavalar, Michael Whitmeyer, and the anonymous referee for helpful comments.

\end{document}